\documentclass[reqno]{article}
\usepackage{amssymb}
\usepackage{latexsym}
\usepackage{amsmath}
\usepackage{amscd}
\newfont{\cyr}{wncyr10}
\newfont{\cyb}{wncyb10}

\begin{document}

\setlength{\textwidth}{8in}%
\setlength{\textheight}{9in}

\def\cc {{\mathfrak c}}
\def\ii {{\mathfrak i}}
\def\UU {{\mathfrak U}}
\def\CC {{\Bbb C}}
\def\HH {{\Bbb H}}
\def\NN {{\Bbb N}}
\def\PP {{\Bbb P}}
\def\QQ {{\Bbb Q}}
\def\RR {{\Bbb R}}
\def\TT {{\Bbb T}}
\def\ZZ {{\Bbb Z}}
\def\sA {{\mathcal A}}
\def\sB {{\mathcal B}}
\def\sC {{\mathcal C}}
\def\sD {{\mathcal D}}
\def\sE {{\mathcal E}}
\def\sF {{\mathcal F}}
\def\sG {{\mathcal G}}
\def\sH {{\mathcal H}}
\def\sI {{\mathcal I}}
\def\sJ {{\mathcal J}}
\def\sK {{\mathcal K}}
\def\sL {{\mathcal L}}
\def\sM {{\mathcal M}}
\def\sN {{\mathcal N}}
\def\sO {{\mathcal O}}
\def\sP {{\mathcal P}}
\def\sQ {{\mathcal Q}}
\def\sR {{\mathcal R}}
\def\sS {{\mathcal S}}
\def\sT {{\mathcal T}}
\def\sU {{\mathcal U}}
\def\sV {{\mathcal V}}
\def\sW {{\mathcal W}}
\def\sX {{\mathcal X}}
\def\sY {{\mathcal Y}}
\def\sZ {{\mathcal Z}}
\def\od {\mathrm{od}}
\def\cf {\mathrm{cf}}
\def\dom {\mathrm{dom}}
\def\id {\mathrm{id}}
\def\int {\mathrm{int}}
\def\cl {\mathrm{cl}}
\def\Hom {\mathrm{Hom}}
\def\ker {\mathrm{ker}}
\def\log {\mathrm{log}}
\def\nwd {\mathrm{nwd}}

\def\TT{{\Bbb T}}
\def\Z{{\Bbb Z}}
\def\ZZ{{\Bbb Z}}

\newcommand{\s }{\mathcal }
\newcommand{\bA}{\mathbf{A}}
\newcommand{\bB}{\mathbf{B}}
\newcommand{\bC}{\mathbf{C}}
\newcommand{\bD}{\mathbf{D}}
\newcommand{\bI}{\mathbf{I}}
\newcommand{\bE}{\mathbf{E}}
\newcommand{\bK}{\mathbf{K}}
\newcommand{\bT}{\mathbf{T}}

\newtheorem{thm}{Theorem}[section]

\newtheorem{theorem}[thm]{Theorem}
\newtheorem{corollary}[thm]{Corollary}
\newtheorem{lemma}[thm]{Lemma}
\newtheorem{claim}[thm]{Claim}
\newtheorem{axiom}[thm]{Axiom}
\newtheorem{conjecture}[thm]{Conjecture}
\newtheorem{fact}[thm]{Fact}
\newtheorem{hypothesis}[thm]{Hypothesis}
\newtheorem{assumption}[thm]{Assumption}
\newtheorem{proposition}[thm]{Proposition}
\newtheorem{criterion}[thm]{Criterion}
\newtheorem{definition}[thm]{Definition}
\newtheorem{definitions}[thm]{Definitions}
\newtheorem{discussion}[thm]{Discussion}
\newtheorem{example}[thm]{Example}
\newtheorem{notation}[thm]{Notation}
\newtheorem{remark}[thm]{Remark}
\newtheorem{remarks}[thm]{Remarks}
\newtheorem{problem}[thm]{Problem}
\newtheorem{terminology}[thm]{Terminology}
\newtheorem{question}[thm]{Question}
\newtheorem{questions}[thm]{Questions}
\newtheorem{notation-definition}[thm]{Notation and Definitions}
\newtheorem{acknowledgement}[thm]{Acknowledgement}


\title{Exactly $n$-resolvable Topological Expansions
}

\author{W.W. Comfort\footnote{Department of Mathematics and Computer Science,
 Wesleyan University, Wesleyan Station,
Middletown, CT 06459; phone: 860-685-2632; FAX 860-685-2571; Email: wcomfort@wesleyan.edu}}
\author{Wanjun Hu\footnote{Department of Mathematics and Computer Science, Albany State University,
Albany, GA 31705; phone: 229-886-4751;Email: Wanjun.Hu@asurams.edu}}

\maketitle



\begin{abstract}
For $\kappa$ a cardinal, a space $X=(X,\sT)$ is $\kappa$-{\it resolvable}
if $X$ admits $\kappa$-many pairwise disjoint $\sT$-dense subsets; 
$(X,\sT)$ is {\it exactly} $\kappa$-{\it resolvable} if it
is $\kappa$-resolvable
but not $\kappa^+$-resolvable.

The present paper complements and supplements the authors' earlier work,
which showed for suitably restricted spaces
$(X,\sT)$ and cardinals $\kappa\geq\lambda\geq\omega$ that $(X,\sT)$, if
$\kappa$-resolvable, admits an expansion $\sU\supseteq\sT$, with
$(X,\sU)$ Tychonoff if $(X,\sT)$ is Tychonoff, such that
$(X,\sU)$ is $\mu$-resolvable for all $\mu<\lambda$ but is not
$\lambda$-resolvable (cf. Theorem~3.3 of \cite{comfhu10}). 
Here the ``finite case" is addressed. The authors show in ZFC for
$1<n<\omega$: (a)~every
$n$-resolvable space $(X,\sT)$ admits an exactly $n$-resolvable
expansion $\sU\supseteq\sT$; (b)~in some cases, even with
$(X,\sT)$ Tychonoff, no choice of $\sU$ is available such
that $(X,\sU)$ is quasi-regular;
(c)~if $n$-resolvable, $(X,\sT)$ admits an exactly
$n$-resolvable quasi-regular expansion $\sU$ if and
only if either $(X,\sT)$ is itself exactly $n$-resolvable and quasi-regular
or $(X,\sT)$ has a subspace which is either $n$-resolvable and
nowhere dense or is $(2n)$-resolvable. In particular, every
$\omega$-resolvable quasi-regular space admits an exactly
$n$-resolvable quasi-regular
expansion.
Further, for many familiar
topological properties $\PP$, one may choose $\sU$ so
that $(X,\sU)\in\PP$ if $(X,\sT)\in\PP$.

\end{abstract}
\maketitle

\section{Introduction}
Let $\kappa>1$ be a (possibly finite) cardinal. Generalizing a concept
introduced by Hewitt~\cite{hewa}, Ceder~\cite{ceder64} defined a
space $(X,\sT)$ to be $\kappa$-{\it resolvable} if there is a family of
$\kappa$-many pairwise disjoint nonempty subsets
of $X$, each $\sT$-dense in $X$.
Generalizations of this concept (for example: the dense sets are perhaps
not pairwise disjoint, but have pairwise intersections which are
``small" in some sense; or, the dense sets are required to be Borel, or
to be otherwise restricted), were introduced and studied
in subsequent decades, for example in \cite{maly74a}
\cite{maly97}, \cite{comfgarc98}, \cite{comfgarc01}, \cite{maly98}.

We refer the reader to such works as \cite{juhss}, \cite{comfhu07},
\cite{juhss2}, \cite{comfhu10} for extensive bibliographic references
relating to the existence of spaces, typically Tychonoff spaces, which
satisfy certain prescribed resolvability properties but not others. The
flavor of our work \cite{comfhu10} is quite different from that of
other papers known to us. In those papers, broadly speaking, the
objective is either (a)~to find conditions on a space sufficient to
ensure some kind of resolvability or (b)~to construct by {\it ad hoc}
means spaces which for certain infinite cardinals $\lambda$ are
$\lambda$-resolvable (sometimes in a modified sense) but which are not
$\kappa$-resolvable for specified $\kappa>\lambda$. In \cite{comfhu10},
in contrast, a broader spectrum of results is enunciated. We showed there that
the tailor-made specific spaces constructed by those {\it ad hoc}
arguments arise as instances of a widely available phenomenon, in this
sense: {\it every} Tychonoff space satisfying mild necessary
conditions admits larger Tychonoff topologies as in (b) above. The
constructions of \cite{comfhu10} are based on the
$\sK\sI\sD$ expansion technique introduced in \cite{huthesis} and
developed further in \cite{hu03}, \cite{comfhu03}, \cite{comfhu04},
\cite{comfhu07}. Roughly speaking, the present work in the finite
context parallels theorems
(cf. \cite{comfhu10}(especially Theorem~3.3)) about $\kappa$-resolvability
when $\kappa$ is infinite. Specifically we show for fixed $n<\omega$
that every
$n$-resolvable space
admits an exactly
$n$-resolvable expansion. In some cases, even when the initial space is
Tychonoff, the expansion cannot be chosen to be quasi-regular.
Further, we characterize explicitly those
$n$-resolvable quasi-regular spaces which do admit an exactly $n$-resolvable
quasi-regular expansion.
It is a pleasing feature of our
arguments that for many familiar topological properties $\PP$,
when the initial hypothesized space $(X,\sT)$
has $\PP$ and does admit an exactly $n$-resolvable quasi-regular expansion
$\sU$, one may arrange also that $(X,\sU)$ has $\PP$.

{\it Ad hoc} constructions of Tychonoff spaces which for fixed
$n<\omega$ are exactly $n$-resolvable have been
available for some time~\cite{vd93b}; see also \cite{cederpear},
\cite{elkin69c}, \cite{eck97}, \cite{fengmasa99}, and \cite{feng00}
for other examples, not all Tychonoff.

\begin{remark}
{\rm
In a preliminary version of this paper, submitted to this journal August
31, 2010, we purported to have proved the statements claimed in 
our abstract \cite{comfhu08abs} and \cite{comfhu10}(5.4(**)). We are
grateful to the referee for indicating a simple counterexample (see
\ref{example} below for a broad generalization of the suggested
argument);
that example helped us to recognize the unavoidable relevance of the
quasi-regularity property which figures prominently in this work, and
to find the more delicate
correct condition captured in Theorem~\ref{lem46} below.
}
\end{remark}

Following van Douwen~\cite{vd93b}, we call a space {\it crowded} if it has
no isolated points. (Some authors prefer the term {\it dense-in-itself}.)
Obviously every resolvable space is crowded. 
\begin{definition}
{\rm
Let $\kappa>1$ be a (possibly finite)
cardinal and let $X=(X,\sT)$ be a space. Then

(a) $X$ is {\it hereditarily $\kappa$-irresolvable} if no nonempty subspace
of $X$ is $\kappa$-resolvable in the
inherited topology;

(b) $X$ is {\it hereditarily irresolvable} if $X$ is hereditarily
$2$-irresolvable; and

(c) $X$
is {\it open-hereditarily irresolvable}
if no nonempty open subspace
of $X$ is resolvable in the
inherited topology.
}
\end{definition}
\begin{notation}
{\rm
(a) Let $(X,\sT)$ be a space and let $Y\subseteq X$. The symbol
$(Y,\sT)$ denotes the set $Y$ with the topology inherited from
$(X,\sT)$.

(b) Given a set $X$ and $\sA\subseteq\sP(X)$, the smallest topology
$\sT$ on $X$ such that $\sT\supseteq\sA$ is denoted $\sT:=\langle\sA\rangle$.
}
\end{notation}
\begin{definition}\label{CCC}
{\rm
Let $\PP$ be a topological property.

(a) $\PP$ is {\it chain-closed} if for each set $X$ and each chain $\sC$ of
$\PP$-topologies on $X$, necessarily $(X,\bigcup\sC)\in\PP$.

(b) $\PP$ is {\it clopen-closed} if:

$[(X,\sT)\in\PP$ and $A\subseteq X$ and
$\sU:=\langle\sT\cup\{A,X\backslash A\}\rangle]
\Rightarrow(X,\sU)\in\PP$.

(c) $\PP$ is {\it $\oplus$-closed} if:

$[(X_0,\sT_0)\in\PP, (X_1,\sT_1)\in\PP, X_0\cap X_1=\emptyset
\Rightarrow (X_0,\sT_0)\oplus(X_1,\sT_1)\in\PP$],\\
\noindent where $(X_0,\sT_0)\oplus(X_1,\sT_1)$ denotes
the ``disjoint union" or ``topological sum"
of the spaces $(X_i,\sT_i)$).

(d) If $\PP$ is a chain-closed and clopen-closed
and $\oplus$-closed property, then $\PP$ is
a {\it CC$\oplus$} property.
}
\end{definition}
\begin{remark}\label{list}
{\rm
We make no attempt to compile a list of all CC$\oplus$ properties
but we note that many familiar topological properties are of that type.
Examples are: $T_0$; $T_1$; $T_2$; quasi-regular; regular; completely regular;
normal; Tychonoff; has a
clopen basis; every two points are separated by a clopen partition; any
concatenation of  CC$\oplus$ properties.

For additional input the interested reader
might consult \cite{engel}(1.5.8).
}
\end{remark}

\begin{discussion}
{\rm
With the necessary preliminaries behind us, we now address the proper
topic of this paper---the search for exactly $n$-resolvable expansions.
This divides naturally and necessarily into two sections:
When the hypothesized topological space is $\omega$-resolvable (``The
Infinite Case"), and when it is not (``The Finite Case"). We treat
these in Sections~2 and 3, respectively.
}
\end{discussion}
\section{The Infinite Case}
We will use frequently the following statement, given by
Illanes~\cite{illanes}(Lemma~2).

\begin{lemma}\label{lem401}
Let $0<n<\omega$. A space which is the union of $n$-many
open-hereditarily irresolvable
subspaces is not $(n+1)$-resolvable.
\end{lemma}

\begin{lemma}\label{lem41}
Let $1<n<\omega$ and let $(X, \sT)$ be an
$\omega$-resolvable space.
Then there exist an expansion $\sT'$ of $\sT$, a nonempty
set $U\in\sT$, and
a $\sT'$-dense partition $\{D_j:j<n\}\cup\{E_j:j<n\}$ of $X$
such that
each $(U\cap D_j, \sT')$ is hereditarily irresolvable.

If in addition $(X,\sT)\in\PP$ with $\PP$ a CC$\oplus$ property, then
$\sT'$ may be
chosen  so that $(X,\sT')\in\PP$.
\end{lemma}
{\sl Proof: }
By transfinite induction we will define an (eventually constant) family
$\{\sT_\eta:\eta<(2^{|X|})^+\}$ of topologies on $X$.

The initial sequence $\{\sT_k:k<\omega\}$ requires special attention.
Recall that the set $\omega$ admits a sequence $\sI_j=\{I^0_j,I^1_j\}$
($j<\omega$)
of two-cell partitions with the property that for each
$F\in[\omega]^{<\omega}$ and $f\in\{0,1\}^F$ one has
$|\bigcap_{j\in F}\,I^{f(j)}_j|=\omega$.
(A quick way to see that is to identify
$\omega$ with a countable dense subset $D$ of the space $\{0,1\}^\omega$
and to set $I^i_j:=D\cap\pi_j^{-1}(\{i\})$ for $j<\omega$,
$i\in\{0,1\}$.) Let $\{S(m):m<\omega\}$ witness the
$\omega$-resolvability
of $(X,\sT)$ and for $j<\omega$, $i\in\{0,1\}$ define

$A^i_j:=\bigcup\{S(m):m\in I^i_j\}$.\\
\noindent Now define $\sT_0:=\sT$ and

$\sT_k:=\langle\sT_0\cup\{A_0^0,A_0^1\}\cup\cdots
\cup\{A^0_{k-1},A^1_{k-1}\}\rangle$
for $1\leq k<\omega$.\\
\noindent Each space $(X,\sT_k)$ is resolvable (in fact,
$\omega$-resolvable) since for $F=\{0,1,\cdots k-1\}$ and
$f\in\{0,1\}^F$ we have $|\bigcap_{j\in F}\,I^{f(j)}_j|=\omega$ so
infinitely many $m$ satisfy $S(m)\subseteq\bigcap_{j\in F}\,A_j^{f(j)}$;
each such $S(m)$ meets each nonempty $U\in\sT=\sT_0$.

Continuing the construction, we define the topologies $\sT_\eta$ for
$\omega\leq\eta<(2^{|X|})^+$. 

For limit ordinals $\eta$, we set
$\sT_\eta:=\bigcup_{\xi<\eta}\,\sT_\xi$.

For successor ordinals $\eta+1$ we have two cases: 
If $(X,\sT_\eta)$ is resolvable we choose a dense partition
$\{A^0_\eta,A^1_\eta\}$ of $(X,\sT_\eta)$ and we set
$\sT_{\eta+1}:=\langle\sT_\eta\cup\{A^0_\eta,A^1_\eta\}\rangle$, and
if $(X,\sT_\eta)$ is irresolvable we set
$\sT_{\eta+1}:=\sT_\eta$.

The definitions of the topologies $\sT_\eta$ are complete.
Routine arguments show that each space $(X,\sT_\eta)$ is crowded, and
Definition~\ref{CCC}(b) and (c) (invoked recursively) shows for
each CC$\oplus$ property $\PP$ that each space
$(X,\sT_\eta)$ has $\PP$ if the initial
space $(X,\sT)$ has $\PP$.

Now for notational simplicity let $\lambda$ be the least ordinal such
that $\sT_\lambda=\sT_{\lambda+1}$ (necessarily with
$\lambda<(2^{|X|})^+$ since for $\eta<\lambda$ we have
$A^0_\eta\in\sT_{\eta+1}\backslash\sT_\eta$). Then
$\lambda\geq\omega$ according to our 
definition of $\{\sT_k: k<\omega\}$.

We set

$R:=\bigcup\{S\subseteq X:(S,\sT_\lambda)$ is resolvable$\}$, and
$W:=X\backslash R$.\\
\noindent Then $W\in\sT_\lambda$, $(W,\sT_\lambda)$ is
hereditarily irresolvable, and $W\neq\emptyset$ since $(X,\sT_\lambda)$
is irresolvable. We fix a nonempty $\sT_\lambda$-basic subset
$U\cap H$ of $W$; here
$U\in\sT=\sT_0$ and
$H=\bigcap_{\eta\in F}\,A_\eta^{f(\eta)}$ for some
$F\in[\lambda+1]^{<\omega}$, $f\in\{0,1\}^F$.

Now let $1<n<\omega$ as hypothesized, choose
$G=\{\eta_j:j<n\}\in[\lambda+1]^n$ such that $G\cap F=\emptyset$, and let
$\{g_j:j<n\}$ be a set of $n$-many distinct functions from $G$ to
$\{0,1\}$. For $j<n$ set $H_j:=\bigcap_{\eta\in G}\,A_\eta^{g_j(\eta)}$
and define

$D_j:=H_j\cap H$ for $j<n$,

$E_j:=H_j\backslash H$ for $1\leq j<n$, and

$E_0:=[H_0\backslash H]\cup[X\backslash
(\bigcup_{j<n}\,D_j\cup\bigcup_{1\leq j<n}\,E_j)]$.\\
\noindent The sets $H$ and $H_j$ ($j<n$) are $\sT_\lambda$-clopen, so
each $D_j$ and $E_j$ ($j<n)$ is $\sT_\lambda$-clopen.

We define

$\sT':=\bigcup\{\sT_\eta:\eta\leq\lambda,\eta\notin F\cup G\}$.

A typical basic open subset of $\sT'$ has the form $U'\cap H'$
with $U'\in\sT=\sT_0$ and $H'=\bigcap _{\eta\in F'}\,A_\eta^{f'(\eta)}$
with $F'\in[\lambda+1]^{<\omega}$, $F'\cap(F\cup G)=\emptyset$ and
$f'\in\{0,1\}^{F'}$; hence the sets
$D_j$, $E_j$ ($j<n$) are dense in $(X,\sT')$. It is clear further that
$\{D_j:j<n\}\cup\{E_j:j<n\}$ is a partition of $X$. It
remains then to show that each space $(U\cap D_j,\sT')$
($j<n$) is hereditarily irresolvable. We note first a weaker statement:
\begin{equation}\label{eq1}
\mbox{Each space~}(U\cap D_j,\sT_\lambda)~(j<n)
\mbox{~is hereditarily irresolvable.}
\end{equation}
\noindent Statement (\ref{eq1}) is clear, since
$U\cap D_j\subseteq U\cap H\subseteq W$ and
$(W,\sT_\lambda)$ is hereditarily irresolvable.

Suppose now for some (fixed) $j<n$ that there is a nonempty set
$A=A^0\cup A^1\subseteq U\cap D_j$ with $\{A^0,A^1\}$
a dense partition of $(A,\sT')$. From (\ref{eq1}), not both $A^0$ and $A^1$ are
dense in $(A,\sT_\lambda)$, so we assume without loss of generality
that $\int_{(A,\sT_\lambda)}\,A^0\neq\emptyset$, say
\begin{equation}\label{eq2}
\emptyset\neq V\cap H'\cap H''\cap A\subseteq A^0
\end{equation}
\noindent with $V\in\sT$, $H'=\bigcap_{\eta\in F'}\,A_\eta^{f'(\eta)}$,
$F'\in[(\lambda+1)\backslash(F\cup G)]^{<\omega}$, $f'\in\{0,1\}^{F'}$,
and with
in addition $H''=\bigcap_{\eta\in F''}\,A_\eta^{f''(\eta)}$, $F''\subseteq F\cup G$,
$f''\in\{0,1\}^{F''}$.

We assume $V\subseteq U$. We assume also, using
$f''|(F\cap F'')=f|(F\cap F'')$ and $f''|(G\cap F'')=g_j|(F\cap F'')$
and replacing
$f''$ by $f\cup g_i$, that $F''=F\cup G$. Then $H''=H_j\cap H=D_j$
and from (\ref{eq2}) we have

$\emptyset\neq V\cap H'\cap H''\cap A=V\cap H'\cap D_j\cap A\subseteq
A^0\subseteq D_j$,\\
\noindent and hence
\begin{equation}\label{eq3}
\emptyset\neq(V\cap H')\cap A\subseteq A^0.
\end{equation}
\noindent Since $V\cap H'\in\sT'$, (\ref{eq3}) shows
$\int_{(A,\sT')}\,A^0\neq\emptyset$, a contradiction since
$A^1=A\backslash A^0$ is dense in $(A,\sT')$.
$\Box$

We use the following lemma only in the case $\kappa=n<\omega$, but
we give the general statement and
proof since these require no additional effort.

\begin{lemma} \label{thm411} Let $\kappa>1$ be a (possibly finite) cardinal
and let $(X, \sT')$ be a space with a dense partition
$\{D_\eta:\eta<\kappa\}\cup\{E_\eta:\eta<\kappa\}$
in which there is a nonempty set $U\in\sT'$ such that
each space $(U\cap D_\eta,\sT')$ is hereditarily irresolvable.
Then there is an expansion $\sT''$ of $\sT'$ such that

{\rm (a)}
$U\cap(\bigcup\{D_\eta:\eta<\kappa\})\in\sT''$; and

{\rm (b)} each set $D_\eta\cup E_\eta$ is dense in $(X,\sT'')$.

If in addition $(X,\sT')\in\PP$ with $\PP$ a CC$\oplus$ property, then
$\sT''$ may be
chosen  so that $(X,\sT'')\in\PP$.
\end{lemma}
{\sl Proof: }
With notation as hypothesized, set
$W:=U\cap(\bigcup\{D_\eta:\eta<\kappa\})$
and define

$\sT'':=\langle\sT'\cup\{W,X\backslash W\}\rangle$.

Clearly (a) holds, since $W=U\cap(\bigcup\{D_\eta:\eta<\kappa\})\in\sT''$.

Definition~\ref{CCC} applies as before to guarantee
that if $(X,\sT')\in\PP$ then also
$(X,\sT'')\in\PP$.

For (b), we fix a nonempty basic set $V''\in\sT''$ and
$\overline{\eta}<\kappa$. We must show
\begin{equation}\label{eq4}
(D_{\overline{\eta}}\cup E_{\overline{\eta}})\cap V''\neq\emptyset.
\end{equation}
For some nonempty set $U'\in\sT''$ we have
either $V''=U'\cap W$ or
$V''=U'\backslash W$.
In the former case since $\emptyset\neq U'\cap U\in\sT'$ and
$D_{\overline{\eta}}$ is
dense in $(X,\sT')$ we have 

$\emptyset\neq(D_{\overline{\eta}}\cap U)\cap
U'=(D_{\overline{\eta}}\cap W)
\cap U'=D_{\overline{\eta}}\cap V''$;\\
\noindent and in the latter case from
$E_{\overline{\eta}}\cap W=\emptyset$
and the density of $E_{\overline{\eta}}$ in $(X,\sT')$
we have $V''\cap E_{\overline{\eta}}=U'\cap
E_{\overline{\eta}}\neq\emptyset$. 
Thus (\ref{eq4}) is proved.
$\Box$

\begin{theorem} \label{cor3-30}
Let $1<n<\omega$ and let $(X, \sT)$ be an
$\omega$-resolvable space.
Then there is an expansion $\sU$ of $\sT$
such that $(X,\sU)$ is exactly $n$-resolvable.

If in addition $(X,\sT)\in\PP$ with $\PP$ a CC$\oplus$ property, then
$\sU$ may be
chosen  so that $(X,\sU)\in\PP$.
\end{theorem}
{\sl Proof: }
Let $\sT'\supseteq\sT$, $U\in\sT$, and $\{D_j:j<n\}\cup\{E_j:j<n\}$ be
as given by Lemma~\ref{lem41}. Then by Lemma~\ref{thm411} there is an
expansion $\sT''$ of $\sT'$ such that

$\emptyset\neq W:=U\cap(\bigcup_{j<n}\{D_j:j<n\})\in\sT''$\\
\noindent and each set $D_j\cup E_j$ is dense in $(X,\sT'')$.

We define $\sU:=\sT''$. As indicated in the statements of
Lemmas~\ref{lem41} and \ref{thm411}, we have $(X,\sU)\in\PP$
if $(X,\sT)\in\PP$.

The family $\{D_j:j<n\}\cup\{E_j:j<n\}$ is a dense partition of
$(X,\sT'')=(X,\sU)$, so $(X,\sU)$ is $n$-resolvable (indeed,
$2n$-resolvable). We have
$W\cap D_j=U\cap D_j\neq\emptyset$ and $(U\cap D_j,\sT')$ is hereditarily
irresolvable,
so (since $\sU\supseteq\sT'$) the space $(W\cap D_j,\sU)$ is hereditarily
irresolvable. The relation $W=\bigcup_{j<n}\,(W\cap D_j)$ expresses $W$ as
the union of $n$-many open-hereditarily irresolvable (even, hereditarily
irresolvable) $\sU$-dense subspaces, so
from Lemma~\ref{lem401} we have that $(W,\sU)$ is not $(n+1)$-resolvable.
Then surely, since $\emptyset\neq W\in\sU$, the space $(X,\sU)$ is not
$(n+1)$-resolvable.
$\Box$
\section{The Finite Case}
We have shown for $1<n<\omega$ that each $\omega$-resolvable space
admits an exactly $n$-resolvable expansion.
That result leaves unresolved the following two questions: (a)~Does
every $n$-resolvable ($\omega$-irresolvable) space admit an exactly
$n$-resolvable expansion? (b)~If not, which $n$-resolvable spaces do
admit such an expansion? In this Section we respond fully to those
questions, as follows.
First, we show in Theorem~\ref{expansionsexist} that every
$(n+1)$-resolvable, $\omega$-irresolvable
space admits an exactly
$n$-resolvable expansion which is not quasi-regular.
Next, we give in Lemma~\ref{neg} a set of
conditions sufficient to ensure that a given $n$-resolvable space admits
no exactly $n$-resolvable quasi-regular expansion. Then, profiting from a referee's
report and leaning heavily on an example given by
Juh\'asz, Soukup and Szentmikl\'ossy~\cite{juhss}, we show in
Theorem~\ref{example} that for every $n>1$ there do exist (many)
Tychonoff spaces satisfying the hypotheses of Lemma~\ref{neg};
thus, not every $n$-resolvable Tychonoff space admits an exactly $n$-resolvable
quasi-regular expansion.
Finally in Theorem~\ref{lem46}, sharpening the results given,
we characterize internally those $n$-resolvable quasi-regular
spaces which do admit an
exactly $n$-resolvable quasi-regular expansion; and we show, as in the
$\omega$-resolvable case treated in Section~2, that
for every CC$\oplus$ property $\PP$ the expansion may be
chosen in $\PP$ if the initial space was in $\PP$.

Twice in this section we will invoke the
following useful result Theorem of Illanes~\cite{illanes}.
We remark
that the anticipated
generalization of Theorem~\ref{illanes}
to (arbitrary) infinite cardinals of countable
cofinality, not needed here,
was given by Bhaskara Rao~\cite{brao94}.
\begin{theorem}\label{illanes}
A space which is $n$-resolvable for each integer $n<\omega$ is
$\omega$-resolvable.
\end{theorem}

We follow Oxtoby~\cite{oxtoby} in adopting the terminology
of this next definition; alternatively,
the spaces in
question might be referred to as spaces with a {\it regular-closed
$\pi$-base}. (We are grateful to Alan Dow and Jerry Vaughan for helpful
correspondence concerning these terms.)

\begin{definition}\label{quasireg}
{\rm
A space $(X,\sT)$ is {\it quasi-regular} if for every nonempty
$U\in\sT$ there
is a nonempty $V\in\sT$ such that
$V\subseteq\overline{V}^{(X,\sT)}\subseteq U$.
}
\end{definition}

The condition of quasi-regularity has the flavor of a (very weak) separation
condition. Clearly every regular space is
quasi-regular. We note that a quasi-regular space need not be a $T_1$-space.

\begin{theorem}\label{expansionsexist}
Let $1<n<\omega$. Every $(n+1)$-resolvable, $\omega$-irresolvable space
admits an exactly $n$-resolvable expansion $(X,\sU)$ which
is not quasi-regular.
\end{theorem}
{\sl Proof: }
It suffices to prove that
there is nonempty $U\in\sT$ such that $(X\backslash U,\sT)$ is
$n$-resolvable
and $(U,\sT)$ admits an exactly $n$-resolvable non-quasi-regular
expansion.
For if $\sU'$ is such a topology on $U$ then
$\sU:=\langle\sT\cup\sU'\cup\{U,X\backslash U\}\rangle$ is as required
for $X$. (Note for clarity: In the notation of Definition~\ref{CCC} we have
$(X,\sU)=(U,\sU')\oplus(X\backslash U,\sT)$.)

Since $(X,\sT)$ is not $\omega$-resolvable, by Theorem~\ref{illanes}
there is $m$ such that $n<m<\omega$ and $(X,\sT)$ is exactly
$m$-resolvable. We set $R:=\bigcup\{S\subseteq X:(S,\sT)$ is
$(m+1)$-resolvable$\}$. Then $(R,\sT)$ is closed in $(X,\sT)$ and
$n$-resolvable (even
$(m+1)$-resolvable), and with $U:=X\backslash R$ we have $\emptyset\neq
U\in\sT$. For clarity, we denote by $\sT'$ the trace of $\sT$ on $U$.
By the previous paragraph it suffices to find an exactly
$n$-resolvable, non-quasi-regular expansion $\sU'$ of $\sT'$ on $U$.

Let $\{D_i:i<m\}$ be a dense partition of $(U,\sT')$, and for $i<m$ set

$R_i:=\bigcup\{S:S\subseteq D_i, (S,\sT')$ is resolvable$\}$.\\
Then $R_i$ is closed in $(D_i,\sT')$, and $(R_i,\sT')$ is resolvable.
We claim for $i<m$ that
\begin{equation}\label{(1)}
R_i \mbox{~is nowhere dense in~} (D_i,\sT').
\end{equation}
Indeed if for some $\overline{i}<m$ there is nonempty $V\in\sT'$ such
that $V\cap D_{\overline{i}}\subseteq R_{\overline{i}}$ then (since
$(R_{\overline{i}},\sT')$ is resolvable and $U\cap R_{\overline{i}}$ is open in
$(R_{\overline{i}},\sT')$) the set $V\cap R_{\overline{i}}=V\cap
D_{\overline{i}}$ is resolvable. Then $V=\bigcup_{i<m}\,(V\cap D_i)$ would
be $(m+1)$-resolvable, contrary to
the fact that $(U,\sT)=(U,\sT')$ is hereditarily $(m+1)$-irresolvable.
Thus (\ref{(1)}) is
proved.

It follows that the set $R(U):=\bigcup_{i<m}\,R_i$ is nowhere dense in
$(U,\sT')$. For $i<m$ we write

$E_i:=D_i\backslash\overline{R(U)}^{(U,\sT')}$.\\
\noindent Then each set $E_i$ is dense in $(U,\sT')$, each space
$(E_i,\sT')$ is hereditarily irresolvable, and
$U=(\bigcup_{i<m}\,E_i)\cup\overline{R(U)}^{(U,\sT')}$.

Now set $E:=\bigcup\{E_i:i<n\}$ and
$\sU':=\langle\{\sT'\cup\{E\}\rangle$.
As the union of $n$-many hereditarily
irresolvable subsets, the space $(E,\sT')$ is not $(n+1)$-resolvable (by
Lemma~\ref{lem401}), hence is exactly $n$-resolvable.
Then since $E$ is open in $(U,\sT')$, also $(U,\sT')$
is exactly $n$-resolvable.

To see that $(U,\sU')$ is not quasi-regular, let $\emptyset\neq
V\in\sU'$ with $V\subseteq E\in\sU'$. There is $W\in\sT'$
(with $W\subseteq U$) such that
$V=W\cap E$, and using $(E,\sU')=(E,\sT')$ we have

$\overline{V}^{(U,\sU')}=\overline{W\cap
E}^{(U,\sU')}=\overline{W}^{(U,\sU')}=\overline{W\cap E}^{(U,\sT')}
=\overline{W}^{(U,\sT')}\supseteq W\cap E_n\neq\emptyset$,\\
\noindent so $\overline{V}^{U,\sU')}\subseteq E$ fails.
$\Box$
We continue with a lemma which is a
routine strengthening of Theorem~\ref{cor3-30}.
\begin{lemma}\label{subspace}
Let $1<n<\omega$ and let $(X,\sT)$ be an $n$-resolvable
space with a nonempty
$\omega$-resolvable subspace. Then there is an expansion $\sU$ of $\sT$
such that $(X,\sU)$ is exactly $n$-resolvable.

If in addition $(X,\sT)\in\PP$ with $\PP$ a CC$\oplus$ property, then
$\sU$ may be
chosen  so that $(X,\sU)\in\PP$.
\end{lemma}
{\sl Proof: }
Let $A\subseteq X$ have the property that $(A,\sT)$ is
$\omega$-resolvable. Replacing $A$ by $\overline{A}^{(X,\sT)}$ if necessary, we
assume that $A$ is $\sT$-closed in $X$. By Theorem~\ref{cor3-30} there
is a topology $\sW$ on $A$, with $(X,\sW)\in\PP$ if $(X,\sT)\in\PP$,
such that $\sW$ expands (the trace of) $\sT$ on $A$ and $(A,\sW)$ is
exactly $n$-resolvable. Both $A$ and $X\backslash A$ are clopen in the topology
$\sU:=\langle\sT\cup\sW\rangle$;
here $(X,\sU)=(A,\sW)\oplus(X\backslash A,\sU)$.
Since $(X\backslash A,\sU)=(X\backslash A,\sT)$ is $n$-resolvable
and $(A,\sU)=(A,\sW)$ is exactly $n$-resolvable,
the space $(X,\sU)$ is exactly $n$-resolvable. When $(X,\sT)\in\PP$ and
$(A,\sU)\in\PP$, necessarily $(X,\sU)\in\PP$
since property $\PP$ is $\oplus$-closed
(see Definition~\ref{CCC}(c).
$\Box$
\begin{definition}\label{-maximal}
{\rm
For $\kappa$ a cardinal, a space $(X, \sT)$ is
$\kappa$-{\it maximal} if no nonempty subspace of $X$
is both nowhere dense and
$\kappa$-resolvable.
}
\end{definition}
\begin{lemma}\label{neg}
Let $1<n<m<2n<\omega$ and let $(X,\sT)$ be an $n$-maximal,
$m$-resolvable space which is hereditarily ($m+1)$-irresolvable.
Then $(X,\sT)$ admits no exactly $n$-resolvable quasi-regular expansion.
\end{lemma}
{\sl Proof:}
Let $\{D_i:i<m\}$ be a dense partition of $(X,\sT)$, and for $i<m$ set

$R_i:=\bigcup\{S:S\subseteq D_i, (S,\sT)$ is resolvable$\}$.\\
Then $R_i$ is closed in $(D_i,\sT)$, and $(R_i,\sT)$ is resolvable.
Exactly as in the proof of (\ref{(1)}) we have
for $i<m$ that

$R_i$ is nowhere dense in $(D_i,\sT)$.

\noindent It follows that the set $R:=\bigcup_{i<m}\,R_i$ is nowhere dense in
$(X,\sT)$. In what follows for $i<m$ we write

$E_i:=D_i\backslash\overline{R}^{(X,\sT)}$.\\
\noindent Then each set $E_i$ is dense in $(X,\sT)$, each space
$(E_i,\sT)$ is hereditarily irresolvable, and

$X=(\bigcup_{i<m}\,E_i)\cup\overline{R}^{(X,\sT)}$
with $(\bigcup_{i<m}\,E_i)\cap\overline{R}^{(X,\sT)}=\emptyset$.

Suppose now that there is an expansion $\sU$ of $\sT$ such that
$(X,\sU)$ is exactly $n$-resolvable and quasi-regular. We claim
\begin{eqnarray}\label{card}
\mbox{there are nonempty~}U''\in\sU\mbox{~and~}F\in[m]^{\leq n}\nonumber\\
\mbox{~such that~}U''\subseteq\bigcup_{i\in F}\,E_i.
\end{eqnarray}
\noindent To prove (\ref{card}) let $R_\sU:=\bigcup\{S\subseteq X:(S,\sU)$ is
$(n+1)$-resolvable$\}$ and set
$V:=X\backslash(R_\sU\cup\overline{R}^{(X,\sT)})$. Then
$\emptyset\neq V\in\sU$ and $(V,\sU)$ is hereditarily $(n+1)$-irresolvable. Now
for $\emptyset\neq U''\subseteq V$ with $U''\in\sU$ we set
$\#(U''):=\{i<m:U\cap E_i\neq\emptyset\}$ and we choose such $U''$ with
$|\#(U'')|$ minimal. Since $U''=\bigcup_{i\in\#(U'')}\,(U''\cap E_i)$ and
$(U'',\sU)$ is $(n+1)$-irresolvable,
if $|\#(U'')|>n$ there is $\overline{i}\in\#(U'')$ such
that $(U''\cap E_{\overline{i}})$ is not dense in $(U'',\sU)$; then some
nonempty $U'\subseteq U''$ with $U'\in\sU$ satisfies $U'\cap
E_{\overline{i}}=\emptyset$, and then $|\#(U')|<|\#(U'')|$. That
contradiction establishes (\ref{card}).

Since $(X,\sU)$ is quasi-regular, there are by (\ref{card}) a set
$F\in[m]^{\leq n}$ and nonempty $U'\in\sU$ such that

$U'\subseteq\overline{U'}^{(X,\sU)}\subseteq\bigcup_{i\in F}\,E_i$.\\
\noindent We choose and fix such $U'$ and we assume further, reindexing
the sets $E_i$ if necessary, that
\begin{equation}\label{U'}
U'\subseteq\overline{U'}^{(X,\sU)}\subseteq\bigcup_{i<n}\,E_i.
\end{equation}
For every nonempty $U\subseteq U'$ such that $U\in\sU$ the space
$(U,\sU)$
is $n$-resolvable, so $(U,\sT)$ is $n$-resolvable; therefore, since
$(X,\sT)$ is $n$-maximal, such $U$ is not nowhere dense in $(X,\sT)$.
Given such $U$, let $V=V(U):=\int_{(X,\sT)}\,\overline{U}^{(X,\sT)}$; then
$\emptyset\neq U\cap V$. Every
nonempty $W\in\sU$ such that $W\subseteq U\cap V$ is $n$-resolvable in
the topology
$\sU$, hence in the topology $\sT$. More explicitly, we claim:
\begin{eqnarray}\label{W}
\noindent\hspace{.3in}U\subseteq U',U\in\sU,V=\int_{(X,\sT)}\,\overline{U}^{(X,\sT)}, W\in\sU,\nonumber\\
\noindent\hspace{.3in}W\subseteq U\cap V,
i<n\Rightarrow W\subseteq\overline{W\cap E_i}^{(X,\sT)}.
\end{eqnarray}
If the claim fails for some such $W$ then for some $i_0<n$ we have, defining
$W':=W\backslash\overline{W\cap E_{i_0}}^{(X,\sT)}$, that $\emptyset\neq
W'\in\sU$ and $W'\cap E_{i_0}=\emptyset$.

Now recursively for $i<n$ we will define a nonempty set $U_i\in\sT$ and
a (possibly empty) set
$F_i\subseteq E_i$ as follows.

For $i=0$, if $W'\cap E_0$ is empty or crowded in the topology $\sT$
we take $U_0=X$ and
$F_0=\emptyset$. If
$W'\cap E_0$ is neither empty nor crowded in the topology $\sT$
we choose a point $x_0$ which
is isolated in $(W'\cap E_0,\sT)$ and we choose $U_0\in\sT$ such that
$U_0\cap W'\cap E_0=\{x_0\}$; then we define $F_0:=\{x_0\}$.

We continue similarly, assuming $i<n$ and that $U_j\in\sT$
and $F_j\subseteq E_j$ have been
defined for all $j<i$ in such a way that $U_{j'}\subseteq U_j$ when
$0\leq j<j'<i$. If $U_{i-1}\cap W'\cap E_i$ is empty or crowded
in the topology $\sT$ we take
$U_i=U_{i-1}$ and $F_i=\emptyset$. If
$U_{i-1}\cap W'\cap E_i$ is neither empty nor crowded in the topology
$\sT$
we choose a point
$x_i$ which
is isolated in $U_{i-1}\cap W'\cap E_i$ and we choose $U_i\in\sT$ such that
$U_i\cap U_{i-1}\cap W'\cap E_i=\{x_i\}$; we assume, replacing $U_i$ by
$U_i\cap U_{i-1}$ if necessary, that $U_i\subseteq U_{i-1}$.
Further in this case we define $F_i:=\{x_i\}$.

The definitions are complete for (all) $i<n$.
We have that $\emptyset\neq U_{n-1}\cap W'\in\sU$, so
$(U_{n-1}\cap W',\sU)$
is $n$-resolvable. It follows then that 
$(U_{n-1}\cap W'\cap E_i,\sT)$ is crowded for at least one $i<n$ (for
otherwise $|U_{n-1}\cap W'\cap E_i|\leq1$ for each
$i<n$, with $U_{n-1}\cap W'\cap E_{i_0}=\emptyset$ and then
$|U_{n-1}\cap W'|\leq n-1$,
so $(U_{n-1}\cap W',\sU)$ cannot be $n$-resolvable).

[We note in passing that no separation properties have been
used or  assumed here, in
particular it is not assumed that a finite subset of $X$ is closed in
$(X,\sU)$. We use only the fact that a space of cardinality $n-1$ or
less  cannot be $n$-resolvable.] 

The relation

$U_{n-1}\cap W'=\bigcup_{i<n}\,(U_{n-1}\cap W'\cap E_i)  $\\
\noindent expresses $U_{n-1}\cap W'$ as the union of $n$-many sets, each
of them hereditarily
irresolvable in the topology $\sT$. (Some of these sets may be
singletons, and at least one of them, namely 
with $i=i_0$ is empty.)
Then by Lemma~\ref{lem401}
the space $(U_{n-1}\cap W',\sT)$ is not $n$-resolvable, so
$(U_{n-1}\cap W',\sU)$ is not $n$-resolvable.
That contradiction
completes the proof of claim (\ref{W}).

Now recursively for $i<n$ we will define nonempty sets
$U_i$, $W_i\in\sU$ and $V_i$, $O_i\in\sT$. To
begin, take $U_0:=U'$, $V_0:=\int_{(X,\sT)}\,\overline{U_0}^{(X,\sT)}$,
and $W_0:=U_0\cap V_0$.
From (\ref{W}) we have
$\overline{W_0}^{(X,\sT)}=\overline{W_0\cap E_0}^{(X,\sT)}$.
Since $(W_0,\sT)$ is $n$-resolvable, also
$(\overline{W_0}^{(X,\sT)},\sT)=(\overline{W_0\cap E_0}^{(X,\sT)},\sT)$
is $n$-resolvable, so
$\int_{(X,\sT)}\,(\overline{W_0\cap E_0}^{(X,\sT)})\neq\emptyset$ (since
$(X,\sT)$ is $n$-maximal). Then since $(E_0,\sT)$ is hereditarily
irresolvable, we have $\int_{(E_0,\sT)}\,(W_0\cap E_0)\neq\emptyset$; we
choose nonempty $O_0\in\sT$ such that
$\emptyset\neq O_0\cap E_0\subseteq W_0\cap E_0$.

We continue similarly, assuming $i<n$ and that
$U_j$, $W_j\in\sU$ and $V_j$, $O_j\in\sT$ have been defined for all $j<i$ so
that $O_{j'}\subseteq O_j$ when $0\leq j<j'<i$. We set $U_i:=W_{i-1}\cap
O_{i-1}$, $V_i:=\int_{(X,\sT)}\,\overline{U_i}^{(X,\sT)}$, and
$W_i:=U_i\cap V_i$. By the preceding
argument (applied now to $W_i\cap E_i$ in place of $W_0\cap E_0$)
there is nonempty $O_i\in\sT$ such that $O_i\cap E_i\subseteq
W_i\cap E_i$. Replacing $O_i$ by $O_i\cap O_{i-1}$ if necessary, we
have $O_i\subseteq O_{i-1}$.

The definitions are complete for (all) $i<n$. We set

$M:=O_{n-1}\backslash(\overline{U'}^{(X,\sU)}\cup\overline{R}^{(X,\sT)})$.\\
\noindent To see that $M\neq\emptyset$, we argue as follows. First,
$\emptyset\neq O_{n-1}\in\sT$ and $\overline{R}^{(X,\sT)}$ is nowhere
dense in $(X,\sT)$, so
$O_{n-1}\backslash\overline{R}^{(X,\sT)}\neq\emptyset$; then, $E_n$ is
dense in $(X,\sT)$, so
$(O_{n-1}\backslash\overline{R}^{(X,\sT)})\cap E_n\neq\emptyset$;
finally, since $\overline{U'}^{(X,\sU)}\cap E_n=\emptyset$, we have

$M\supseteq(M\cap E_n)=((O_{n-1}\backslash\overline{R}^{(X,\sT)})
\backslash\overline{U'}^{(X,\sU)})\cap E_n\neq\emptyset$.

For $i<n$ we have $M\cap(O_{n-1}\cap E_i)=\emptyset$, so
\begin{equation}\label{M}
M=\bigcup_{n\leq i<m}\,(M\cap E_i).
\end{equation}
\noindent  Since $m<2n$, relation (\ref{M}) expresses $M$ as the union of
fewer than $n$-many subsets, each hereditarily irresolvable in the
topology $\sT$, so again by Lemma~\ref{lem401}
the space $(M,\sT)$ is not $n$-resolvable. Then also $(M,\sU)$ is not
$n$-resolvable, which with the relation $\emptyset\neq M\in\sU$
contradicts the assumption that $(X,\sU)$ is (exactly) $n$-resolvable.
$\Box$

The following argument shows, as promised, that, for every integer $n>1$,
Tychonoff spaces satisfying the conditions of Lemma~\ref{neg} exist in
profusion. 
\begin{theorem}\label{example}
Let $1<n<\omega$ and let $\kappa\geq\omega$. Then there is an
$n$-resolvable Tychonoff space $(X,\sT)$ such that $|X|=\kappa$
and $(X,\sT)$ admits no exactly
$n$-resolvable quasi-regular expansion.
\end{theorem}

{\sl Proof: }
According to \cite{juhss}(4.1) there is a dense subspace $Y$ of the
space $\{0,1\}^{2^\kappa}$ such that $|Y|=\kappa$ and $Y$ is {\it
submaximal} in the sense that every dense subspace of $Y$ is open in
$Y$. Let $\{D_i:i<n+1\}$ be a family of pairwise disjoint dense
subspaces of $\{0,1\}^{2^\kappa}$, each homeomorphic to $Y$,
and set $X:=\bigcup_{i<n+1}\,D_i$. (To
find such spaces $D_i$, let $H$ be the subgroup of $\{0,1\}^{2^\kappa}$
generated by $Y$, note that $|H|=|Y|=\kappa$, choose any $(n+1)$-many
cosets of H and chose a dense homeomorph of $Y$ inside each of those.)

Clearly $X=(X,\sT)$ is $n$-resolvable, indeed $(n+1)$-resolvable. Hence 
to see that $(X,\sT)$ admits no exactly $n$-resolvable
quasi-regular expansion it
suffices, according to Lemma~\ref{neg} (taking $m=n+1$ there), to show that 
$(X,\sT)$ is $n$-maximal and hereditarily $(n+2)$-irresolvable.

In fact, $(X,\sT)$ is even $2$-maximal. To see that, let $A$ be nowhere
dense in $(X,\sT)$. Then $A\cap D_i$ is nowhere dense in $(D_i,\sT)$ for
each $i<n+1$, hence is hereditarily closed, hence is (hereditarily)
discrete. The relation $A=\bigcup_{i<n+1}\,(A\cap D_i)$ expresses $A$ as
the union of finitely many discrete sets, so (since $(X,\sT)$ is a
$T_1$-space) the space $(A,\sT)$ is not crowded and hence is not
$2$-resolvable.

To see that $(X,\sT)$ is hereditarily $(n+2)$-irresolvable, we begin with
a preliminary observation.
\begin{equation}\label{submax}
\mbox{Every submaximal space is hereditarily irresolvable}.
\end{equation}
For the proof, let $S\subseteq Z$ with $Z$ a submaximal space, and suppose
that $\{S_0,S_1\}$ is a dense partition of $S$. Then $S':=Z\backslash
S_0=S_1\cup(Z\backslash S)$ is dense in $Z$ and hence open in $Z$, so
$S'\cap S=S_1$ is open in $S$; then $S_0=S\backslash S_1$ cannot be
dense in $S$, a contradiction completing the proof of (\ref{submax}).

Suppose now in the present case that $A\subseteq X$ and that $(A,\sT)$ is
$(n+2)$-resolvable. Replacing $A$ by $\overline{A}^{(X,\sT)}$
if necessary, we assume
that $A$ is closed in $(X,\sT)$. We consider two cases.

\underline{Case 1}. $A$ is nowhere dense in $(X,\sT)$. Then $A\cap D_i$ is
nowhere dense in $(D_i,\sT)$ for each $i<n+1$, hence is
discrete and closed in $D_i$. Then $A=\bigcup_{i<n+1}\,(A\cap D_i)$
is not crowded, hence is
not resolvable.

\underline{Case 2}.  Case 1 fails. Then there is nonempty $U\in\sT$ such
that $U\subseteq A$. Each set $U\cap D_i$ with $i<n+1$ is dense in
$(U,\sT)$, and each space $(U\cap D_i,\sT)$ is hereditarily irresolvable
since (by (\ref{submax})) each space $(D_i,\sT)$ is hereditarily irresolvable.
The relation $U=\bigcup_{i<n+1}\,(U\cap D_i)$ then
expresses the $(n+2)$-resolvable space
$(U,\sT)$ as the union of $(n+1)$-many dense, hereditarily
irresolvable subsets. That contradicts Lemma~\ref{lem401}.
$\Box$

Finally in Theorem~\ref{lem46} we give the promised internal
characterization of those $n$-resolvable quasi-regular
spaces which admit an exactly
$n$-resolvable quasi-regular expansion.

\begin{theorem}\label{lem46}
Let $1<n<\omega$ and let $(X,\sT)$ be an $n$-resolvable
quasi-regular space. Then
these conditions are equivalent.

{\rm (a)} $(X,\sT)$ admits an exactly $n$-resolvable quasi-regular expansion;

{\rm (b)} either
\begin{itemize}
\item[(i)] $(X,\sT)$ is exactly $n$-resolvable; or
\item[(ii)] $(X,\sT)$ is not $n$-maximal; or
\item[(iii)] $(X,\sT)$ is not hereditarily $(2n)$-irresolvable.
\end{itemize}

If in addition conditions {\rm (a)} and {\rm (b)} are satisfied and
$(X,\sT)\in\PP$ with $\PP$ a CC$\oplus$ property, then
the exactly $n$-resolvable quasi-regular expansion $\sU$ of $\sT$ may be
chosen  so that $(X,\sU)\in\PP$.
\end{theorem}
{\sl Proof: }

(b)$\Rightarrow$(a). If (i) holds there is  nothing to prove. We\
assume that (i) fails and we show
(ii)$\Rightarrow$(a) and (iii)$\Rightarrow$(a).

Let $\{D_i:i<n\}$ be a dense partition of
$(X,\sT)$ and let $(A,\sT)$ be a subspace of $(X,\sT)$ such that
either $(A,\sT)$ is $n$-resolvable and nowhere dense in
$(X,\sT)$, or $(A,\sT)$ is ($2n$)-resolvable. Replacing $A$ by
$\overline{A}^{(X,\sT)}$ if necessary, we assume
that $A$ is $\sT$-closed in $X$. If $(A,\sT)$ is $\omega$-resolvable the
desired conclusion is given by Lemma~{\ref{subspace}}
(using the fact that quasi-regularity is a CC$\oplus$ property),
so we assume that
$(A,\sT)$ is not $\omega$-resolvable; then by Lemma~\ref{illanes} there is
$m$ such that $n< m<\omega$ and $(A,\sT)$ is exactly
$m$-resolvable. Let $\{E_i:i<m\}$ witness that fact, and set
$E:=\bigcup_{i<n}\,E_i$. We claim

\begin{equation}\label{8}
\mbox{If~} \sU \mbox{~is an expansion of~} \sT \mbox{~in which~} E
\mbox{~is~} \sU\mbox{-clopen,
then~} (X,\sU) \mbox{~is not~} (n+1)\mbox{-resolvable.}
\end{equation}

Indeed, if $(X,\sU)$ is $(n+1)$-resolvable then its clopen subset
$(E,\sU)$ admits a dense partition of the form
$\{F_j:j<n+1\}$; then $\{F_j:j<n+1\}\cup\{E_i:n\leq i<m\}$ would be a dense
$(m+1)$-partition of $(A,\sT)$, contrary to the fact that $(A,\sT)$ is
not $(m+1)$-resolvable. Thus (\ref{8}) is proved.

Suppose now that $(A,\sT)$ is nowhere dense in $(X,\sT)$, as in (ii),
and set $\sU:=\langle\sT\cup\{E,X\backslash E\}\rangle$. Each
nonempty set
$U\in\sU$ meets either $E$ or $X\backslash E$, hence
(since $\int_{(X,\sT)}\,A=\emptyset$)
meets either $E$ or $X\backslash A$; hence $U$ 
meets each set $D_i$ ($i<n$) or each set $E_i$ ($i<n$),
so $\{D_i\cup E_i:i<n\}$
witnesses the fact that $(X,\sU)$ is $n$-resolvable. It then follows from
(\ref{8}), since $E$ is $\sU$-clopen, that $(X,\sU)$ is exactly $n$-resolvable.

Suppose that $(A,\sT)$ is
($2n$)-resolvable, as in (iii), set
$\sU:=\{\sT\cup\{X\backslash A,E,A\backslash E\}\rangle$, and
let $\emptyset\neq U\in\sU$. Clearly
$U$ meets either
$X\backslash A$ or $E$ or $A\backslash E$. It follows, since
$(X\backslash A,\sU)=(X\backslash A,\sT)$ and $(E,\sU)=(E,\sT)$ and
$(A\backslash E,\sU)=(A\backslash E,\sT)$, that $U$ meets either each set
$D_i$ ($i<n$) or each set $E_i$ ($i<n$) or each set $E_{n+i}$ ($i<n$).
Thus $\{D_i\cup E_i\cup E_{n+i}:i<n\}$
witnesses the fact that $(X,\sU)$ is $n$-resolvable. It then follows from
(\ref{8}) that $(X,\sU)$ is exactly $n$-resolvable.

Every CC$\oplus$ property $\PP$
is preserved under passage from $(X,\sT)$ to $(X,\sU)$ under the
constructions given in (ii) and (iii), so in particular $(X,\sU)$ is
quasi-regular since $(X,\sT)$ was assumed quasi-regular. 

(a)~$\Rightarrow$~(b).
Assume that (b) fails. For $i<n$ let

$R_i:=\bigcup\{S\subseteq X:(S,\sT)$ is
$(n+i+1$)-resolvable$\}$,\\
\noindent and note that

$R_0\supseteq R_1\supseteq\ldots\supseteq R_{i-1}\supseteq R_i
\ldots\supseteq R_{n-1}$

\noindent with $R_0=X$ (since $(X,\sT)$ is $(n+1)$-resolvable)
and $R_{n-1}=\emptyset$ (since $(X,\sT)$ is hereditarily
$(2n)$-irresolvable).

Since each set $R_i$ ($i<n-1$) is closed in $(X,\sT)$,
each set $R_i\backslash R_{i+1}$ ($i<n-1$) is open in $(R_i,\sT)$;
further, $R_i\backslash R_{i+1}$ is
$(n+i+1)$-resolvable and hereditarily $(n+i+2)$-irresolvable.
For $i<n-1$ let
$U_i:=\int_{(X,\sT)}\,R_i$, and $B_i:=R_i\backslash U_i$.

Let $\sU$ be a refinement of $\sT$ such that $(X,\sU)$
is exactly $n$-resolvable, and set

$R:=\bigcup\{S\subseteq X:S$ is $(n+1)$-resolvable$\}$.\\
\noindent Then with $W:=X\backslash R$ we have $\emptyset\neq W\in\sU$,
and $(W,\sU)$ is $n$-resolvable and hereditarily $(n+1)$-irresolvable.

Since
$X=\bigcup_{i<n-1}\,(R_i\backslash R_{i+1})$, we have
$W=\bigcup_{i<n-1}\,(W\cap(R_i\backslash R_{i+1}))$. 
Each $B_i$ ($i<n-1$) is closed and nowhere dense in $(R_i,\sT)$, hence
in $(X,\sT)$, so $B:=\bigcup_{i<n-1}\,B_i$ and each of its subsets
is nowhere dense in $(X,\sT)$; since $(X,\sT)$ is $n$-maximal, $(B,\sT)$
and each of its subsets is $n$-irresolvable. Then since $(W,\sT)$ is
$n$-resolvable the relation $W\subseteq B$ fails, so there is $i<n-1$
such that $V:=W\cap(U_i\backslash R_{i+1})$ is nonempty.

From $i<n-1$ it follows with $m:=n+i+1$ that

$1<n<m=n+i+1<n+n-1+1=2n$\\
\noindent and the proof concludes with the observation that
the following two facts
contradict (the instance $X=V$ of)
Lemma~\ref{neg}.

(A) the space $(V,\sT)$ is $n$-maximal, $(n+i+1)$-resolvable, and
hereditarily $(n+i+2)$-irresolvable; and

(B) $\emptyset\neq V\in\sU$, $(V,\sU)$ is $n$-resolvable,
quasi-regular, and
not $(n+1)$-resolvable; that is, the space $(V,\sU)$
is exactly $n$-resolvable and quasi-resolvable. 

$\Box$

\end{document}